\documentclass{amsart}
\usepackage{amsmath}
\usepackage{amsthm}
\usepackage{graphicx} 
\usepackage[backend=biber, style=numeric]{biblatex}
\usepackage[a4paper,margin=1in]{geometry} 
\newcommand{\stirling}[2]{%
  \genfrac{[}{]}{0pt}{}{#1}{#2}%
}

\newtheorem{definition}{Definition}

\newtheorem*{lemma}{Lemma}
\newtheorem*{proposition}{Proposition}
\newtheorem*{corollary}{Corollary}
\newtheorem{theorem}{Theorem}

\title[Cycles in random generalized derangements]{The probability that two random generalized derangements have the same number of cycles}
\author{Dmitrii Smorchkov \& Mikl\'os B\'ona}

\begin{document}

\begin{abstract}
    Using singularity analysis, we show that the probability that two independently and uniformly selected generalized derangements of size $n$ have the same number of cycles is asymptotically the same as the corresponding probability for unrestricted permutations, namely, $$\sim \frac{1}{2\sqrt{\pi\log n}}\;.$$
\end{abstract}
\maketitle

\section{Introduction}
When studying large combinatorial objects, a natural question arises: if we select two random objects of size n from a given class, what is the probability that they share a specified shape, profile or parameter? For example, the authors of \cite{MR2737320} study the probability that two random compositions have the same number of parts, while Wilf \cite{wilf2005variancestirlingcyclenumbers} studies the corresponding question for the number of cycles in random permutations.

Recall that the Stirling numbers of the first kind $\stirling{n}{k}$ count  permutations of size $n$ that have exactly $k$ cycles. In \cite{wilf2005variancestirlingcyclenumbers} Wilf showed that for two independently and randomly selected permutations the probability that they have the same number of cycles is asymptotically
$$\frac{1}{(n!)^2}\sum\limits_{k=1}^n \stirling{n}{k}^2\sim \frac{1}{2\sqrt{\pi\log n}}\;.$$
Wilf’s elegant argument relies on the polynomial structure of the generating function for the Stirling numbers and on specific properties of the Gamma function. This approach does not extend directly to permutations with restrictions on cycle lengths.
\begin{definition}
Let $S$ be a (possibly empty) finite set of positive integers. An {\em $S$-derangement} is a permutation whose cycle lengths do \textbf{not} belong to $S$. 
\end{definition}

We denote the collection of all $S$-derangements  of length $n$ by
$\mathcal{D}_n^S$.

Using this notation, $\mathcal{D}_n^\emptyset$ is the set of all permutations of size $n$, and $\mathcal{D}_n^{\{1\}}$ is the set of ordinary derangements, that is, permutations without fixed points. To simplify the notation, we fix $S$ throughout the remainder of this paper and omit the superscript.

In this paper, we prove the perhaps surprising fact that the probability that two randomly selected generalized $S$-derangements of size $n$ have the same number of cycles is asymptotically the same as the corresponding probability for unrestricted permutations. This is interesting, because the possible numbers of cycles are much less numerous. For instance, an ordinary derangement
of length $n$ can have at most $\lfloor n/2 \rfloor$ cycles, whereas for $S=\{1,2,\cdots, k\}$, the highest possible number of cycles is $\lfloor n/(k+1) \rfloor$. 

In order to keep the paper self-contained while not slowing down the flow of the proof too much, out of the three theorems we use, we include two in the main part of the article, and we present the third in the Appendix. 

\section{The proof of our main result}

\begin{definition}
We refine $\mathcal{D}_n$ according to the number of cycles by setting
$$\mathcal{D}_{n,k}=\{\,p \in \mathcal{D}_n\,\vert\,\ \text{\#cycles}(p)=k\,\}\,.$$
Write $D_n=|\mathcal{D}_n|$ and $D_{n,k}=|\mathcal{D}_{n,k}|$ for the corresponding cardinalities.
\end{definition}
Our goal is to show that $$\frac{1}{(D_n)^2}\sum\limits_{k=1}^n(D_{n,k})^2 \sim \frac{1}{2\sqrt{\pi\log n}}$$
for any finite $S$.
We use the methods of \cite{Flajolet_Sedgewick_2009} to establish a Gaussian local limit law for the appropriately normalized cycle counts as n tends to infinity.
\begin{definition}
   Let $F(z,u)$ be  the bivariate generating function given by $$F(z, u)=\sum\limits_{n=0}^\infty\sum\limits_{k=0}^nD_{n,k}u^k\frac{z^n}{n!}$$
    with the convention that $D_{0,0}=1$ and $D_{n,0}=0$ for all positive integers $n$. We denote the horizontal layers of $F(z,u)$  by $$f_n(u)=[z^n]F(z,u)=\frac{1}{n!}\sum\limits_{k=0}^nD_{n,k}u^k\,.$$
\end{definition}

\begin{definition}
    For $n$ such that $D_n>0$, let $X_n$ be the number of cycles in a random $S$-derangement of size $n$. Its probability generating function is $$p_n(u)=\mathbb{E}u^{X_n}=\sum\limits_{k=0}^n\frac{D_{n,k}}{D_n}u^k=\frac{f_n(u)}{f_n(1)}\;.$$
\end{definition}
Using the symbolic properties of generating functions from \cite{bóna2006walk} we derive the explicit formula for $F(z,u)$ as a \textit{collection of cycles marked by parameter $u$}. Thus, $$F(z,u)=\exp\left(\sum\limits_{j\notin S}u\frac{z^j}{j} \right)=\exp{\left(u\log\frac{1}{1-z}-u\sum\limits_{j\in S}\frac{z^j}{j}\right)}=\exp{\left(-u\sum\limits_{j\in S}\frac{z^j}{j}\right)}(1-z)^{-u}\,.$$
The series $\sum_{j\notin S}z^j/j$ defines a holomorphic function $g(z)$ on the open unit disk. Consequently, $F(z,u) = \exp (ug(z))$ is holomorphic on the set $\{z\in \mathbb{C}\;\vert\;|z| < 1\}\times\mathbb{C}$.

Following the proof of the singularity analysis theorem for variable-exponent perturbations \cite[Theorem IX.11]{Flajolet_Sedgewick_2009} (see our Appendix for a detailed discussion of the theorem and how we apply it here), we obtain
\begin{equation} \label{perturb} f_n(u)=[z^n]F(z,u)=\frac{\exp{\left(-u\sum\limits_{j\in S}\frac{1}{j}\right)}}{\Gamma(u)}n^{u-1}\left(1+O(n^{-\frac{1}{2}})\right) \end{equation}
uniformly for $u$ in a sufficiently small neighborhood $\Omega$ of $1$. For example, $\Omega=\{u\in\mathbb{C}\;\vert\;|u-1|\leq1/3\}$ will work. Consequently, uniformly on $\Omega$, we have $$p_n(u)=\frac{[z^n]F(z,u)}{[z^n]F(z, 1)}=\frac{\exp{\left((1-u)\sum\limits_{j\in S}\frac{1}{j}\right)}}{\Gamma(u)}n^{u-1}\left(1+O(n^{-\frac{1}{2}})\right)\;.$$

We first recall the Quasi-powers Theorem
\cite[Theorem IX.8]{Flajolet_Sedgewick_2009}.
For a function $f$ analytic at $u=1$, with $f(1)\neq0$, define
its analytic mean and variance by
\[
m(f)=\frac{f'(1)}{f(1)},
\qquad
v(f)=\frac{f''(1)}{f(1)}
+\frac{f'(1)}{f(1)}
-\left(\frac{f'(1)}{f(1)}\right)^2.
\]

\begin{theorem}[Quasi-powers theorem]
Let $X_n$ be a sequence of non-negative discrete random variables
with probability generating functions $P_n(u)$. Suppose that,
uniformly in a fixed complex neighborhood of $u=1$,
\[
P_n(u)=A(u)\cdot B(u)^{\beta_n}
\left(1+O\left(\frac{1}{\kappa_n}\right)\right),
\]
where $\beta_n,\kappa_n\to+\infty$, the functions $A(u)$ and
$B(u)$ are analytic at $u=1$, and $A(1)=B(1)=1$.
Suppose also that the variability condition
\[
v(B)=B''(1)+B'(1)-(B'(1))^2\neq0
\]
holds. Then the mean and variance of $X_n$ satisfy
\begin{align*}
\mu_n&=\mathbb{E}X_n
=\beta_n m(B)+m(A)+O(\kappa_n^{-1}),\\
\sigma_n^2&=\mathbb{V}X_n
=\beta_n v(B)+v(A)+O(\kappa_n^{-1}).
\end{align*}
Moreover, the distribution of $X_n$, after standardization,
is asymptotically Gaussian:
\[
\mathbb{P}\left\{\frac{X_n-\mu_n}{\sigma_n}\leq x\right\}
=\Phi(x)+O\left(\kappa_n^{-1}+\beta_n^{-1/2}\right),
\]
uniformly for $x\in\mathbb{R}$, where $\Phi$ is the distribution
function of the standard normal distribution.
\end{theorem}

We now apply this theorem to the number $X_n$ of cycles in a
random $S$-derangement of size $n$.
Writing $n^{u-1}=\left(e^{u-1}\right)^{\log n}$, the uniform
estimate for $p_n(u)$ obtained above has the required
quasi-power form, with
\begin{align*}
\beta_n&=\log n,\\
\kappa_n&=\sqrt{n},\\
A(u)&=\frac{\exp{\left((1-u)\sum\limits_{j\in S}
\frac{1}{j}\right)}}{\Gamma(u)},\\
B(u)&=\exp{(u-1)}.
\end{align*}
Both sequences $\beta_n$ and $\kappa_n$  tend to infinity, and both functions $A(u)$ and $B(u)$ are analytic
in a neighborhood of $u=1$, with $A(1)=B(1)=1$.
Furthermore,
\[
B^{\prime\prime}(1)+B^{\prime}(1)-(B^\prime(1))^2=1\neq0,
\]
so the variability condition is satisfied.
Thus, all the hypotheses of the theorem hold.

Using the properties of derivatives of the Gamma function and substituting these functions and sequences into the theorem's
mean and variance formulas gives
\begin{align*}
\mu_n&=\mathbb{E}X_n
=\log n+\gamma-\sum\limits_{j\in S}\frac{1}{j}
+O(\kappa_n^{-1}),\\
\sigma_n^2&=\mathbb{V}X_n
=\log n+\gamma-\frac{\pi^2}{6}
-\sum\limits_{j\in S}\frac{1}{j}
+O(\kappa_n^{-1}).
\end{align*}
Consequently,
\[
\frac{X_n-\mu_n}{\sigma_n}
\xrightarrow[n\to\infty]{d}\mathcal{N}(0,1),
\]
with an error term
\[
O(\kappa_n^{-1}+\beta_n^{-1/2})
=O((\log n)^{-1/2}).
\]


This Gaussian limit is not new: V. Goncharov \cite{MR10922} established it for unrestricted permutations in the 1940s by analyzing characteristic functions. The convergence to a Gaussian distribution for generalized derangements is stated in \cite{Flajolet_Sedgewick_2009}. We nevertheless verify the hypotheses of the Quasi-powers theorem because they are also required by the local limit theorem \cite[Theorem IX.14]{Flajolet_Sedgewick_2009}. We include this theorem below for easy reference.

\begin{theorem}[Local Limit Law] \label{lll}
Let $X_n$ be a sequence of non-negative discrete random variables
with probability generating functions $p_n(u)$. Suppose that
$p_n(u)$ satisfies the hypotheses of the Quasi-powers Theorem.
In particular, the approximation
\[
p_n(u)=A(u)\cdot B(u)^{\beta_n}
\left(1+O\left(\frac{1}{\kappa_n}\right)\right)
\]
holds uniformly in a sufficiently small fixed complex neighborhood $\Omega$ of $1$.
Suppose, in addition, that there exists $K>1$ such that
\[
|p_n(u)|\leq K^{-\beta_n}
\]
uniformly for all $u$ in the intersection of the unit circle
and the complement $\mathbb{C}\setminus\Omega$.
Then the distribution of $X_n$ satisfies a local limit law
of the Gaussian type, with a rate of convergence
$O(\beta_n^{-1/2}+\kappa_n^{-1})$.
More precisely, writing
$\mu_n=\mathbb{E}X_n$,
$\sigma_n^2=\mathbb{V}X_n$, and
$\phi(x)=(2\pi)^{-1/2}e^{-x^2/2}$, we have
\[
\sup\limits_{x\in\mathbb{R}}
\left|
\sigma_n\mathbb{P}\left\{
X_n=\lfloor\mu_n+x\sigma_n\rfloor
\right\}
-\phi(x)
\right|
=O(\beta_n^{-1/2}+\kappa_n^{-1}).
\]
\end{theorem}

To apply this theorem, we must also estimate $p_n(u)$ on the unit circle, rather than only in a small neighborhood of $1$.

\begin{lemma}
    For $n$ sufficiently large, for each $0<\delta<\pi$, there exist $C>0$ and $K>1$ such that $$|p_n(u)|<CK^{-\log n}$$ for all $u$ satisfying $|u|=1$ and $|\arg u|>\delta$.
\end{lemma}
\begin{proof} Recall that $p_n(u)=f_n(u)/f_n(1)$, where $f_n(u)=[z^n]F(z,u)$. As $F(z,u)$ is the product of two factors, we will estimate those two factors individually. 
We write $F(z,u)=G(z,u)H(z,u)$ where
\begin{align*}
    G(z,u)=(1-z)^{-u}=\sum\limits_{n=0}^\infty g_n(u)z^n \\
    H(z,u)=\exp{\left(-u\sum\limits_{j\in S}\frac{z^j}{j}\right)} =\sum\limits_{n=0}^\infty h_n(u)z^n\\
    g_n(u)=\frac{u(u+1)\cdots(u+n-1)}{n!} \\
    h_n(u)=\frac{1}{n!}\frac{\partial^n}{\partial z^n}H(z,u)\Big\vert_{z=0}
\end{align*}

Setting $u=e^{it}$ with $-\pi<t\leq \pi$, we obtain the following expression for the squared magnitude of $g_n(u)$:
\begin{equation*}
    |g_n(e^{it})|^2=\prod\limits_{k=1}^{n-1}\frac{|e^{it}+k|^2}{(k+1)^2}=\prod\limits_{k=1}^{n-1}\frac{k^2+2k\cos{t}+1}{(k+1)^2}= \prod\limits_{k=1}^{n-1}\left(1-\frac{2k(1-\cos{t})}{(k+1)^2}\right).
\end{equation*}

Put $\alpha=1-\cos t\in[0,2]$. Since
$0\leq 2k\alpha/(k+1)^2\leq1$ for every $k\geq1$, we may apply
$1-x\leq e^{-x}$ for $0\leq x\leq1$ directly to the factors in the product. This gives
\[
    |g_n(e^{it})|^2
    \leq\exp\left(-2\alpha\sum_{k=1}^{n-1}\frac{k}{(k+1)^2}\right).
\]
Moreover,
\[
    \sum_{k=1}^{n-1}\frac{k}{(k+1)^2}
    =\sum_{j=1}^{n}\left(\frac1j-\frac1{j^2}\right)
    \geq\log n-\frac{\pi^2}{6},
\]
since $\sum_{j=1}^{n}1/j\geq\log n$ and
$\sum_{j=1}^{n}1/j^2\leq\pi^2/6$. Taking square roots, we obtain
\[
    |g_n(e^{it})|
    \leq e^{\alpha\pi^2/6}n^{-\alpha}
    \leq C_2n^{\cos t-1},
\]
where $C_2=e^{\pi^2/3}$. This bound holds uniformly for
$-\pi<t\leq\pi$, including $t=\pi$.

To bound $h_n(u)$, note that $H(z,u)$ is an entire mapping $\mathbb{C}^2\to\mathbb{C}$ and for each fixed $u$ the function $z\mapsto H(z,u)$ is entire. Fix $R>1$. Setting $M_R=\max\limits_{|z|\leq R\;|u|=1}|H(z,u)|$, we obtain the following bound from Cauchy's integral formula:
$$|h_n(u)|=\left|\frac{1}{2\pi i}\int\limits_{|z|=R}\frac{H(z,u)}{z^{n+1}}{\rm d}z\right|\leq M_RR^{-n}\;.$$
This bound holds for all $u$ with $|u|=1$.

Finally, we show that convolving a polynomially decaying sequence with an exponentially decaying sequence preserves polynomial decay. Indeed, 
\begin{align*}
  |f_n(u)|
  &=\left|\sum\limits_{k=0}^ng_k(u)h_{n-k}(u)\right| \\
  &\leq\sum\limits_{k=0}^{\lfloor{n/2}\rfloor}|g_k(u)||h_{n-k}(u)|
  +\sum\limits_{k=\lfloor n/2\rfloor+1}^n|g_k(u)||h_{n-k}(u)| \\
  &\leq \sum\limits_{k=0}^{\lfloor{n/2}\rfloor}1\cdot M_RR^{k-n}
  +\sum\limits_{k=\lfloor n/2\rfloor+1}^nC_2\left(\frac{n}{2}\right)^{\cos{t}-1}\cdot M_RR^{k-n} \\
  &\leq M_R R^{-n}\sum\limits_{k=0}^{\lfloor n/2\rfloor}R^k
  +C_2M_R\left(\frac{n}{2}\right)^{\cos{t}-1}\sum\limits_{k=0}^{\lfloor n/2\rfloor}R^{-k} \\
  &\leq M_RR^{-n}\frac{R^{(n/2)+1}-1}{R-1}
  +C_2M_R\frac{2^{1-\cos{t}}}{1-R^{-1}}n^{\cos{t}-1} \\
  &\leq C_3R^{-n/2}+ C_4\left(e^{1-\cos{t}}\right)^{-\log n}.
\end{align*}
The second term asymptotically dominates the first.

Since $1-\cos{t}>1-\cos{\delta}>0$ and $f_n(1)\to\exp(-\sum_{j\in S}1/j)>0$ as $n\to\infty$ by setting $u=1$ in (\ref{perturb}), we may take $K=e^{1-\cos{\delta}}$ and choose $C>0$ such that $$|p_n(e^{it})|=\frac{|f_n(e^{it})|}{f_n(1)}<CK^{-\log n}\;.$$
\end{proof}

We have already verified the hypotheses of the Quasi-powers
theorem for the number $X_n$ of cycles in a random
$S$-derangement, with
$\beta_n=\log n$ and $\kappa_n=\sqrt{n}$.
In order to apply the Local Limit Law, it remains to verify the additional bound on the unit circle.

Choose $0<\delta<\pi$ sufficiently small that
\[
\{e^{it}\;|\;|t|\leq\delta\}\subseteq\Omega.
\]
The preceding lemma gives constants $C>0$ and $K>1$ such that
\[
|p_n(u)|<CK^{-\log n}
\]
for all $u$ on the unit circle with $|\arg u|>\delta$.
In particular, this bound holds on the part of the unit circle
outside $\Omega$.
Choosing any $\Tilde{K}$ with $1<\Tilde{K}<K$, we have, for all sufficiently
large $n$,
\[
CK^{-\log n}\leq \Tilde{K}^{-\log n}=\Tilde{K}^{-\beta_n}.
\]
Thus, the additional hypothesis of the Local Limit Law
is satisfied.

Applying the Local Limit Law yields the following estimate, which is
the key ingredient in determining the probability that two
random $S$-derangements have the same number of cycles.

\begin{proposition}
    Let $p_{n,k}=\mathbb{P}\{X_n=k\}=\frac{D_{n,k}}{D_n}$ and let $\phi(x)$ denote the standard Gaussian density function, that is, $\phi(x)=(2\pi)^{-1/2}e^{-x^2/2}$. Then
    $$\sup\limits_{x\in \mathbb{R}}|\sigma_n \mathbb{P}\{X_n=\lfloor\mu_n+x\sigma_n\rfloor\}-\phi(x)|=O\left(\frac{1}{\sqrt{\log n}}\right)\,.$$
    In particular, $$p_{n,k}=\frac{1}{\sigma_n}\left(\phi\left(\frac{k-\mu_n}{\sigma_n}\right)+\epsilon_{n,k}\right)$$ where $\epsilon_n=\max_{k\in[n]}|\epsilon_{n,k}|=O\left(\frac{1}{\sqrt{\log n}}\right)$.
\end{proposition}

We now prove the main result of this paper.
\begin{theorem}
    The equality
$$\lim\limits_{n\to\infty}\sigma_n\sum\limits_{k=1}^np_{n,k}^2=\int\limits_{-\infty}^{\infty}\phi^2(x){\rm d}x=\frac{1}{2\sqrt{\pi}}$$
holds.
\end{theorem}
\begin{proof}
   The crucial idea is to view the expression on the left-hand side as a sequence of Riemann sums converging to the integral. To make this argument rigorous, we first truncate the tails of the distribution of $X_n$ and establish convergence to the integral over the interval $[-A,A]$ for fixed A. We then let $A$ tend to infinity.
   
    Fix $A>0$. Define the set of indices $$I_n^A=\{k\in[n]\;|\;\frac{|k-\mu_n|}{\sigma_n}\leq A\}\,.$$ Note that $2\sigma_nA - 1<|I^A_n|\leq2\sigma_nA+1$ if $n$ is large enough. Writing $\phi_{nk}=\phi\left(\frac{k-\mu_n}{\sigma_n}\right)$, we obtain
    \begin{align*}
        \sigma_n\sum\limits_{k\in I^A_n}p_{nk}^2=\sigma_n\sum\limits_{k\in I^A_n}\sigma_n^{-2}(\phi_{nk}+\epsilon_{nk})^2=\sigma_n^{-1}\sum\limits_{k\in I^A_n}\phi_{nk}^2+\sigma_n^{-1}\sum\limits_{k\in I^A_n}(2\phi_{nk}\epsilon_{nk}+\epsilon_{nk}^2) \\ =\sigma_n^{-1}\sum\limits_{k\in I^A_n}\phi_{nk}^2 +o(1)
    \end{align*}
    because
    $$\left|\sigma_n^{-1}\sum\limits_{k\in I^A_n}(2\phi_{nk}\epsilon_{nk}+\epsilon_{nk}^2) \right| \leq \frac{2\sigma_nA+1}{\sigma_n}\left(\frac{2}{\sqrt{2\pi}}\epsilon_n+\epsilon_n^2\right)\to0$$ as $n\to\infty$ for any fixed $A$.

    Meanwhile, $\sigma^{-1}_n\sum_{k\in I_n^A}\phi_{nk}^2$ is a Riemann sum for the continuous function $\phi^2(x)$ on the interval $[-A,A]$ with mesh size $\sigma_n^{-1}\to0$ as $n\to\infty$. Therefore,
    $$\sigma_n\sum\limits_{k\in I^A_n}p_{nk}^2\to\int\limits_{-A}^A\phi^2(x){\rm d}x$$ as $n\to\infty$.

   To bound the tails, note that $\max_{k\in[n]}p_{nk}<\frac{C}{\sigma_n}$ for some $C>0$, by the Proposition. Then we apply Chebyshev’s inequality:
    \begin{align*}
        \sigma_n\sum\limits_{k\notin I^A_n}p_{nk}^2\leq\sigma_n\cdot\max_{k\in[n]}p_{nk}\cdot \sum\limits_{k\notin I^A_n}p_{nk}\leq C\cdot\mathbb{P}\left\{\frac{|X_n-\mu_n|}{\sigma_n}>A\right\}\leq\frac{C}{A^2}\;.
    \end{align*}
    We now obtain lower and upper bounds for the full sum. The lower bound $$\sigma_n\sum\limits_{k=1}^np_{nk}^2\geq\sigma_n\sum\limits_{k\in I^A_n}p_{nk}^2=\int\limits_{-A}^A\phi^2(x){\rm d}x+o(1)$$ implies 
    $$\liminf_{n\to\infty}\sigma_n\sum\limits_{k=1}^np_{nk}^2\geq\int\limits_{-A}^A\phi^2(x){\rm d}x\;.$$
    Similarly, the upper bound
$$\sigma_n\sum\limits_{k=1}^np_{nk}^2=\sigma_n\sum\limits_{k\in I^A_n}p_{nk}^2+\sigma_n\sum\limits_{k\notin I^A_n}p_{nk}^2\leq\int\limits_{-A}^A\phi^2(x){\rm d}x+o(1)+\frac{C}{A^2}$$
    implies
    $$\limsup_{n\to\infty} \sigma_n \sum_{k=1}^n p_{nk}^2\leq\int\limits_{-A}^A\phi^2(x){\rm d}x+\frac{C}{A^2}$$
    for every fixed $A>0$.
    Letting $A\to\infty$ and noting that $\phi^2(x)$ is finitely integrable over the real line we get:
    \begin{align*} \liminf_{n\to\infty}\sigma_n\sum\limits_{k=1}^np_{nk}^2=\limsup_{n\to\infty}\sigma_n\sum\limits_{k=1}^np_{nk}^2=\lim_{n\to\infty}\sigma_n\sum\limits_{k=1}^np_{nk}^2=\int\limits_{-\infty}^\infty\phi^2(x){\rm d}x \\=\frac{1}{2\pi}\int\limits_{-\infty}^\infty e^{-x^2}{\rm d}x=\frac{1}{2\sqrt{\pi}}\;
    \end{align*}
    proving our claim.
\end{proof}
Since $\sigma_n\sim \sqrt{\log n}$ we conclude that
$$\sum\limits_{k=1}^np_{nk}^2\sim\frac{1}{2\sqrt{\pi\log n}}\;.$$

\section*{A generalization}
For the main result we utilized the fact that $\sigma_np_{nk}$ is approximated by $\phi_{nk}$ with the uniformly controlled error and that $\phi^2(x)\in L^1(\mathbb{R})$. But higher powers of $\phi(x)$ are also integrable. Thus, the proof can be naturally extended for the scenario of multiple permutations. Fix an integer $d\geq2$. Given a tuple of $d$ generalized $S$-derangements $(p_1,p_2,\dots,p_d)$ selected uniformly at random with replacement\footnote{The same asymptotic formula holds for sampling without replacement. See the Appendix for a proof.} from $\mathcal{D}_n$ we show that the probability that they all have the same number of cycles satisfies the following asymptotics:
$$\sum\limits_{k=1}^np_{nk}^d\sim\frac{1}{\sqrt{d(2\pi\log n)^{d-1}}}\;.$$

To show this we prove the following theorem.
\begin{theorem}
    Let $d\geq2$ be a fixed integer. The equality
$$\lim\limits_{n\to\infty}\frac{1}{\sigma_n}\sum\limits_{k=1}^n(p_{nk}\sigma_n)^d=\int\limits_{-\infty}^{\infty}\phi^d(x){\rm d}x=\frac{1}{\sqrt{d(2\pi)^{d-1}}}$$
holds.
\end{theorem}
\begin{proof}
   The argument is analogous to that of the proof for $d=2$. Using the same notation, we highlight the inequalities that are different. First, note that
    \begin{align*}
        \sigma_n^{-1}\sum\limits_{k\in I^A_n}(\sigma_np_{nk})^d=\sigma_n^{-1}\sum\limits_{k\in I^A_n}\phi_{nk}^d+\sigma_n^{-1}\sum\limits_{k\in I^A_n}\sum\limits_{r=1}^d \binom{d}{r}\phi_{nk}^{d-r}\epsilon_{nk}^r \\ =\sigma_n^{-1}\sum\limits_{k\in I^A_n}\phi_{nk}^d +o(1)
    \end{align*}
  because $0\leq\phi_{nk}\leq 1$ implies that, for all integers
$a\geq 0$ and $b\geq 1$,
\[
\phi_{nk}^{a}\left|\epsilon_{nk}^{b}\right|
\leq \epsilon_n^{b}\longrightarrow 0
\]
uniformly in $k$, and $|I_n^A|\leq 2\sigma_n A+1$. For the tail estimate we have
    \begin{align*}
         \sigma_n^{-1}\sum\limits_{k\notin I^A_n}(\sigma_np_{nk})^d\leq(\sigma_n\max_kp_{nk})^{d-1}\sum\limits_{k\notin I^A_n}p_{nk}\leq\frac{C^{d-1}}{A^2}\;.
    \end{align*}
    Hence we obtain identical bounds for $\liminf$ and $\limsup$ yielding the limit result that
    was to be proved.
\end{proof}
As an application of the generalized result we can consider $S=\emptyset$. Then having $p_{nk}=\frac{1}{n!}\stirling{n}{k}$, we immediately obtain the following asymptotics for the sum of $d$-th powers of Stirling numbers.
\begin{corollary}
    As $n\to\infty$, the asymptotic equivalence
    $$\sum\limits_{k=1}^n\stirling{n}{k}^d\sim\frac{(n!)^d}{\sqrt{d(2\pi\log n)^{d-1}}}\;$$
    holds.
\end{corollary}

\section*{Appendix}

We present the variable-exponent theorem of Flajolet and Sedgewick
\cite[Theorem~IX.11, pp.~669--670]{Flajolet_Sedgewick_2009} and explain
how its proof gives the coefficient estimate used in Section~2.
We write $\mathcal A(z,u)$ and $\mathcal B(z,u)$ for the analytic
functions in that theorem, to distinguish them from the functions
$A(u)$ and $B(u)$ in the Quasi-powers theorem.

\subsection{The theorem}

\noindent\textbf{Theorem IX.11 (Variable exponent perturbation).}
Let $F(z,u)$ be analytic at $(0,0)$, with nonnegative Taylor
coefficients. Suppose that it admits a representation
\[
F(z,u)=\mathcal A(z,u)+\mathcal B(z,u)C(z)^{-\alpha(u)}
\]
with the following properties, for some $r>\rho>0$ and $\varepsilon>0$:
\begin{enumerate}
\item The functions $\mathcal A$ and $\mathcal B$ are analytic on a
neighborhood of $\{(z,u):|z|\le r,\ |u-1|\le\varepsilon\}$.
\item The function $\alpha$ is analytic for $|u-1|\le\varepsilon$,
and $\alpha(1)\in\mathbb R\setminus\{0,-1,-2,\ldots\}$.
\item The function $C$ is analytic for $|z|\le r$. Its only zero
in this disk is the simple zero $\rho$, and
$\mathcal B(\rho,1)\ne0$.
\item The variability condition $\alpha'(1)+\alpha''(1)\ne0$ holds.
\end{enumerate}
The representation is understood by analytic continuation, with a
consistent choice of the power of $C(z)$. If $X_n$ has probability
generating function
\[
p_n(u)=\frac{[z^n]F(z,u)}{[z^n]F(z,1)},
\]
then, writing $\mu_n=\mathbb E X_n$ and
$\sigma_n^2=\operatorname{Var}(X_n)$, we have
\begin{align*}
\mu_n&=\alpha'(1)\log n+O(1),\\
\sigma_n^2&=\bigl(\alpha'(1)+\alpha''(1)\bigr)\log n+O(1),
\end{align*}
and
\[
\mathbb P\!\left\{\frac{X_n-\mu_n}{\sigma_n}\le x\right\}
=\Phi(x)+O\bigl((\log n)^{-1/2}\bigr)
\]
uniformly for $x\in\mathbb R$, where $\Phi$ is the standard normal
distribution function.

\subsection{The coefficient estimate in the proof}

Put $d=-\rho C'(\rho)$. Since $\rho$ is a simple zero, the local
expansions at $z=\rho$ are
\[
C(z)=d(1-z/\rho)\bigl(1+O(1-z/\rho)\bigr),
\qquad
\mathcal B(z,u)=\mathcal B(\rho,u)+O(1-z/\rho),
\]
uniformly for $u$ near $1$. Thus the leading singular term of $F$ is
\[
\mathcal B(\rho,u)d^{-\alpha(u)}(1-z/\rho)^{-\alpha(u)}.
\]
We estimate the  coefficients of $F(z,u)$  using
\[
[z^n](1-z/\rho)^{-\alpha(u)}
=\rho^{-n}\frac{\Gamma(n+\alpha(u))}
{\Gamma(\alpha(u))\Gamma(n+1)}
=\rho^{-n}\frac{n^{\alpha(u)-1}}{\Gamma(\alpha(u))}
\bigl(1+O(n^{-1})\bigr).
\]
The contribution of $\mathcal A$ is exponentially smaller, while
the remaining singular terms contain an additional factor
$1-z/\rho$. Uniform singularity analysis, after restricting the
neighborhood of $u=1$ if necessary, therefore gives the estimate
established in the proof of Theorem~IX.11:
\begin{equation}\label{eq:app-general-coefficient}
[z^n]F(z,u)=c(u)\rho^{-n}n^{\alpha(u)-1}
\bigl(1+O(n^{-1/2})\bigr),
\qquad
c(u)=\frac{\mathcal B(\rho,u)d^{-\alpha(u)}}
{\Gamma(\alpha(u))}.
\end{equation}
Here $c$ is analytic and nonzero on a sufficiently small
neighborhood of $1$, and the error is uniform there.
Dividing \eqref{eq:app-general-coefficient} by its value at $u=1$
yields
\[
p_n(u)=\frac{c(u)}{c(1)}
\exp\!\bigl((\alpha(u)-\alpha(1))\log n\bigr)
\bigl(1+O(n^{-1/2})\bigr).
\]
This is the quasi-powers approximation from which the Gaussian
conclusion follows. The analytic mean and variance of its
factor $\exp(\alpha(u)-\alpha(1))$ are $\alpha'(1)$ and
$\alpha'(1)+\alpha''(1)$, respectively, which explains the
leading terms in the theorem's mean and variance formulas.
In our paper, it is the intermediate
coefficient estimate \eqref{eq:app-general-coefficient} that is
used first.

\subsection{Application to generalized derangements}

Set
\[
P(z)=\sum_{j\in S}\frac{z^j}{j},
\qquad s=P(1)=\sum_{j\in S}\frac1j,
\qquad H(z,u)=e^{-uP(z)}.
\]
Our generating function is $F(z,u)=H(z,u)(1-z)^{-u}$, so the
theorem applies with
\[
\mathcal A(z,u)=0,\qquad \mathcal B(z,u)=H(z,u),\qquad
C(z)=1-z,\qquad \alpha(u)=u,\qquad \rho=1.
\]
Indeed, $P$ is a polynomial because $S$ is finite, so $H$ is
entire. The function $C$ has the unique simple zero $1$, while
$H(1,1)=e^{-s}>0$, and
$\alpha(1)=1$ and $\alpha'(1)+\alpha''(1)=1$.
The nonnegativity of the coefficients of $F$ follows from its
counting interpretation, or directly from
\[
F(z,u)=\exp\!\left(u\sum_{\substack{j\ge1\\j\notin S}}
\frac{z^j}{j}\right).
\]
Also $d=-\rho C'(\rho)=1$, so $c(u)=e^{-us}/\Gamma(u)$.
Consequently, \eqref{eq:app-general-coefficient} becomes
\begin{equation}\label{eq:app-coefficient}
f_n(u)=[z^n]F(z,u)
=\frac{e^{-us}}{\Gamma(u)}n^{u-1}
\bigl(1+O(n^{-1/2})\bigr).
\end{equation}
This is precisely the coefficient estimate used in (\ref{perturb}).

For completeness, the uniformity on the particular disk
$\Omega=\{u:|u-1|\le1/3\}$ used there can be seen directly from
the same argument. Since $H$ is entire, write
\[
H(z,u)=H(1,u)+(1-z)J(z,u),
\]
where the apparent singularity of $J$ at $z=1$ is removable.
Thus $J$ is entire, and
\[
F(z,u)=e^{-us}(1-z)^{-u}+(1-z)^{1-u}J(z,u).
\]
Choose the branch of $\log(1-z)$ on
$\mathbb C\setminus[1,\infty)$. The terms are analytic in a
common $\Delta$-domain at $1$. The uniform transfer argument
used in the proof of Theorem~IX.11 gives
\[
[z^n]\bigl((1-z)^{1-u}J(z,u)\bigr)
=O\bigl(n^{\operatorname{Re}u-2}\bigr),
\qquad u\in\Omega.
\]
The constants can be chosen uniformly because $J$ is bounded
on the relevant compact sets and both the real and imaginary
parts of $u$ range over bounded intervals. For the leading
term, uniformly on $\Omega$,
\[
[z^n](1-z)^{-u}
=\frac{\Gamma(n+u)}{\Gamma(u)\Gamma(n+1)}
=\frac{n^{u-1}}{\Gamma(u)}\bigl(1+O(n^{-1})\bigr).
\]
Since $e^{-us}/\Gamma(u)$ is bounded away from zero on $\Omega$,
these estimates give a uniform relative error $O(n^{-1})$,
which in particular implies the $O(n^{-1/2})$ error retained
in \eqref{eq:app-coefficient}.

Finally, setting $u=1$ in \eqref{eq:app-coefficient} gives
\[
f_n(1)=e^{-s}\bigl(1+O(n^{-1/2})\bigr)\longrightarrow e^{-s}>0,
\]
because $\Gamma(1)=1$ and $n^{1-1}=1$. Dividing by this
normalizing factor, uniformly for $u\in\Omega$, we obtain
\[
p_n(u)=\frac{f_n(u)}{f_n(1)}
=\frac{e^{(1-u)s}}{\Gamma(u)}n^{u-1}
\bigl(1+O(n^{-1/2})\bigr).
\]
Writing $n^{u-1}=(e^{u-1})^{\log n}$ gives exactly the
quasi-powers parameters used in the paper:
\[
A(u)=\frac{e^{(1-u)s}}{\Gamma(u)},\qquad
B(u)=e^{u-1},\qquad
\beta_n=\log n,\qquad \kappa_n=\sqrt n.
\]
In particular, $A(1)=B(1)=1$ and
$B''(1)+B'(1)-(B'(1))^2=1\neq 0$, as required.
\subsection{Sampling without replacement.}
 The same asymptotic result holds when $(p_1,p_2,\dots,p_d)$ is sampled without replacement. Indeed, sampling without replacement is equivalent to sampling with replacement conditioned on the event that the $p_i$ are all distinct. Define the events $Q_{ij}=\{p_i\neq p_j\}$ for $1\leq i < j \leq d$ and let $Q = \cap_{i<j}Q_{ij}$. Let $C=\{\#cycles(p_1)=\#cycles(p_2)=\dots=\#cycles(p_d)\}$. By the law of total probability,
$$\mathbb{P}(C)=\mathbb{P}(C\,\vert\,Q)\mathbb{P}(Q)+\mathbb{P}(C\,\vert\,Q^c)\mathbb{P}(Q^c)=\mathbb{P}(C\,\vert\,Q)(1-\mathbb{P}(Q^c))+\mathbb{P}(C\,\vert\,Q^c)\mathbb{P}(Q^c)\;.$$
Rearranging terms, we get
\begin{align}\label{diff-prob}
    |\mathbb{P}(C)-\mathbb{P}(C\,\vert\, Q)|=\mathbb{P}(Q^c)|\mathbb{P}(C\,\vert\,Q^c)-\mathbb{P}(C\,\vert\,Q)|\leq\mathbb{P}(Q^c)\;.
\end{align}
We now show that $\mathbb{P}(Q^c)$ is asymptotically negligible relative to $\mathbb{P}(C)$. By the union bound,
$$\mathbb{P}(Q^c)=\mathbb{P}\left(\cup Q_{ij}^c\right)\leq\sum\limits_{1\leq i < j \leq d}\mathbb{P}(Q^c_{ij})=\binom{d}{2}\frac{1}{D_n}\;.$$
Recall that $D_n/n!=[z^n]F(z,1)=[z^n]\left(\frac{\exp(-P(z))}{1-z}\right)$ and that $F(z,1)$ has a unique singularity, a simple pole at $z=1$. Applying singularity analysis, we conclude that $D_n/n!\to \exp(-P(1))=e^{-s}$ as $n\to\infty$. Dividing both sides of \eqref{diff-prob} by $\mathbb{P}(C)$, we obtain 
$$\left|1-\frac{\mathbb{P}(C\,|\,Q)}{\mathbb{P}(C)}\right|\leq\frac{\mathbb{P}(Q^c)}{\mathbb{P}(C)}=O\left(\frac{(\log n)^{(d-1)/2}}{n!}\right)\to 0$$
as $n\to\infty$. Thus, $\mathbb{P}(C\,|\,Q)\sim\mathbb{P}(C)$.
\printbibliography
\end{document}